\documentclass[11pt]{article}
\usepackage{amscd}
\usepackage{amsfonts}
\usepackage{amsmath}
\usepackage{amssymb}
\usepackage{amsthm}
\usepackage{bbm}
\usepackage{CJK}
\usepackage{fancyhdr}
\usepackage{graphicx}
\usepackage{hyperref}
\usepackage{indentfirst}
\usepackage{latexsym}
\usepackage{mathrsfs}
\usepackage{xypic}

\usepackage[top=1in,bottom=1in,left=1.25in,right=1.25in]{geometry}
\def\<{\langle}
\def\>{\rangle}

\date{}
\allowdisplaybreaks

\begin{document}

\renewcommand{\baselinestretch}{1.2}
\renewcommand{\arraystretch}{1.0}

\title{\bf Heisenberg double as braided commutative Yetter-Drinfel'd module algebra over Drinfel'd double in multiplier Hopf algebra case}
\author
{
  \textbf{Tao Yang} \footnote{Corresponding author. College of Science, Nanjing Agricultural University, Nanjing 210095, Jiangsu, CHINA.
             E-mail: tao.yang.seu@gmail.com} \ ,
  \textbf{Xuan Zhou} \footnote{Department of Mathematics, Jiangsu Institute of Education, Nanjing 210013, Jiangsu, CHINA} \ ,
  \textbf{Juzhen Chen} \footnote{College of Science, Nanjing Forestry University, Nanjing 210037, Jiangsu, CHINA}
}
\maketitle

\begin{center}
\begin{minipage}{12.cm}
 {\bf Abstract}
 Based on a pairing of two regular multiplier Hopf algebras $A$ and $B$, Heisenberg double $\mathscr{H}$ is the smash product $A \# B$
 with respect to the left regular action of $B$ on $A$.
 Let $\mathscr{D}=A\bowtie B$ be the Drinfel'd double, then Heisenberg double $\mathscr{H}$ is a Yetter-Drinfel'd $\mathscr{D}$-module algebra,
 and it is also braided commutative by the braiding of Yetter-Drinfel'd module,
 which generalizes the results in \cite{S11} to some infinite dimensional cases.
\\

\\


\end{minipage}
\end{center}
\normalsize

\section{Introduction}
\def\theequation{\thesection.\arabic{equation}}
\setcounter{equation}{0}

 As shown in \cite{S11}, for a finite dimensional Hopf algebra $H$, there is a Yetter-Drinfel'd $\mathcal{D}(H)$-module algebra structure
 on the Heisenberg double $\mathcal{H}(H^{*})$ endowed with a heterotic action of the Drinfel'd double $\mathcal{D}(H)$,
 moreover $\mathcal{H}(H^{*})$ is braided commutative in terms of the braiding of Yetter-Drinfel'd module.

 One question naturally arises: if this result also holds for some infinite dimensional Hopf algebras?

 As we know, the duality of an infinite dimensional Hopf algebra is no longer a Hopf algebra.
 Therefore, the classic Hopf theory obviously cannot deal with this question,
 we should use another method: multiplier Hopf algebra theory.
 Multiplier Hopf algebras, considered as a generalization of Hopf algebras,
 play a very important role in the duality of a class of infinite dimensional Hopf algebras (see \cite{V94, V98}).
 Many constructions, such as Drinfel'd double (see \cite{DV04,YW11a}), Yetter-Drinfel'd module (see \cite{D08,YW11b}), have naturally generalized.
 Multiplier Hopf algebra becomes a useful tool to deal with some infinite dimensional Hopf algebra questions.

 In this paper, we use the multiplier Hopf algebra theory to deal with a more general case, and give a positive answer to the question.
 Based on a pairing of regular multiplier Hopf algebras, we show the Heisenberg double
 $\mathscr{H}$ is a (left-left) $\mathscr{D}$-Yetter-Drinfel'd module algebra, and braided commutative.
 Furthermore, we apply the conclusion to some special infinite dimensional Hopf algebras, such as co-Frobenius Hopf algebras.

 The paper is organized in the following way.
 In section 2, we recall some notions which we will use in the following, such as
 multiplier Hopf algebras, Drinfel'd doubles and Yetter-Drinfel'd modules.

 In section 3, we define an action and a coaction of Drinfel'd double $\mathscr{D}$ on Heisenberg double $\mathscr{H}$, which makes
 $\mathscr{H}$ a (left-left) $\mathscr{D}$-Yetter-Drinfel'd module algebra (see Theorem 3.4).
 Moreover, the Heisenberg doube $\mathscr{H}$ is a braided $\mathscr{D}$-commutative algebra.
 And $\mathscr{H}$ is the braided product $\mathscr{H} = A \propto B$, where $A$ and $B$ are braided commutative
 Yetter-Drinfel'd $\mathscr{D}$-module algebras by restriction.

 In section 4, we apply the results as above to the usual Hopf algebras and derive some interesting results,
 including the main results in paper \cite{S11} as a corollary.

\section{Preliminaries}
\def\theequation{\thesection.\arabic{equation}}
\setcounter{equation}{0}

 We begin this section with a short introduction to multiplier Hopf algebras.

 Throughout this paper, all spaces we considered are over a fixed field $K$ (such as field $\mathbb{C}$ of complex numbers).
 \emph{Algebras may or may not have units}, but always should be non-degenerate, i.e., the multiplication maps (viewed as bilinear forms) are non-degenerate.
 For an algebra $A$, the multiplier algebra $M(A)$ of $A$ is defined as the largest algebra with unit in which $A$ is a dense ideal
 (see the appendix in \cite{V94}).

 Now, we recall the definitions of a multiplier Hopf algebra (see \cite{V94} for details).
 A comultiplication on algebra $A$ is a homomorphism $\Delta: A \longrightarrow M(A \otimes A)$ such that $\Delta(a)(1 \otimes b)$ and
 $(a \otimes 1)\Delta(b)$ belong to $A\otimes A$ for all $a, b \in A$. We require $\Delta$ to be coassociative in the sense that
 \begin{eqnarray*}
 (a\otimes 1\otimes 1)(\Delta \otimes \iota)(\Delta(b)(1\otimes c))
 = (\iota \otimes \Delta)((a \otimes 1)\Delta(b))(1\otimes 1\otimes c)
 \end{eqnarray*}
 for all $a, b, c \in A$ (where $\iota$ denotes the identity map).

 A pair $(A, \Delta)$ of an algebra $A$ with non-degenerate product and a comultiplication $\Delta$ on $A$ is called a multiplier Hopf algebra,
 if the linear maps $T_{1}, T_{2}: A\otimes A \longrightarrow A\otimes A$ defined by
 $T_{1}(a\otimes b)=\Delta(a)(1 \otimes b)$ and $T_{2}(a\otimes b)=(a \otimes 1)\Delta(b)$
 are bijective.

 A multiplier Hopf algebra $(A, \Delta)$ is called regular if $(A, \Delta^{cop})$ is also a multiplier Hopf algebra,
 where $\Delta^{cop}$ denotes the co-opposite comultiplication defined as $\Delta^{cop}=\tau \circ \Delta$ with $\tau$ the usual flip map
 from $A\otimes A$ to itself (and extended to $M(A\otimes A)$). In this case, $\Delta(a)(b \otimes 1), (1 \otimes a)\Delta(b) \in A \otimes A$
 for all $a, b\in A$.
 By Proposition 2.9 in \cite{V98}, multiplier Hopf algebra $(A, \Delta)$ is regular if and only if the antipode $S$ is bijective from $A$ to $A$.

 If $(A, \Delta)$ is a multiplier Hopf algebra and if $A$ has an identity, then by Theorem 2.4 in \cite{V98} $(A, \Delta)$ is a Hopf algebra.
 We will use the adapted Sweedler notation for regular multiplier Hopf algebras (see \cite{V08}).
 We will e.g., write $\sum a_{(1)} \otimes a_{(2)}b$ for $\Delta(a)(1 \otimes b)$ and $\sum ab_{(1)} \otimes b_{(2)}$ for $(a \otimes 1)\Delta(b)$,
 sometimes we omit the $\sum$.

 \subsection{Pairing and Drinfel'd double}

 Start with two regular multilier Hopf algebras $(A, \Delta)$ and $(B, \Delta)$ together with a non-degenerate bilinear map
 $\langle \cdot, \cdot \rangle$ from $A\times B$ to $K$.
 By Definition 2.8 in \cite{DrV01}, this bilinear map is called a pairing if certain conditions are fulfilled.
 The main property is that the product in $A$ is dual to the coproduct in $B$ and vice versa.
 There are however certain regularity conditions, needed to give a correct meaning of this statement.
 The investigation of these conditions was done in \cite{DrV01}.

 For $a\in A, b\in B$, then recall from Section 1.2 in \cite{D04}, $a\blacktriangleright b$, $b\blacktriangleright a$,
 $a\blacktriangleleft b$ and $b\blacktriangleleft a$ can be defined in the following way.
 Take $a'\in A$ and $b'\in B$, the left multiplications are defined by the formulas
 \begin{eqnarray*}
 (b\blacktriangleright a)a' = \<a_{(2)}, b \> a_{(1)}a', &&  (a\blacktriangleright b)b' = \<a, b_{(2)}\> b_{(1)}b', \\
 (a\blacktriangleleft b)a' = \<a_{(1)}, b \> a_{(2)}a', &&  (b\blacktriangleleft a)b' = \<a, b_{(1)}\> b_{(2)}b'.
 \end{eqnarray*}
 These formulas make sense because $A$ is a regular multiplier Hopf algebra.
 The right multiplications are defined similarly.

 The regularity conditions (six equivalent conditions listed in Proposition 2.7 of \cite{DrV01})
 on the dual paring say that the multipliers $b\blacktriangleright a$ and $a\blacktriangleleft b$ in $M(A)$
 (resp. $a\blacktriangleright b$ and $b\blacktriangleleft a$ in $M(B)$) actually belong to $A$ (resp. $B$).
 Then as shown in section 2.1 of \cite{DV04}, it is possible to state that the product and the coproduct are dual to each other:
 \begin{eqnarray*}
 \<a, bb'\> = \<b'\blacktriangleright a, b\> = \<a\blacktriangleleft b, b'\>, \\
 \<aa', b\> = \< a, a'\blacktriangleright b\> = \<a', b\blacktriangleleft a\>.
 \end{eqnarray*}
 In this way we get four modules and all these modules are unital.

 A stronger result however is possible here, coming from the existence of local units, see e.g. Proposition 2.2 in \cite{DrVZ99}.
 For instance, take $b\in B$, then there are elements $\{a_{1}, a_{2},\ldots, a_{n}\}$ in $A$ and $\{b_{1}, b_{2},\ldots, b_{n}\}$ in $B$
 such that $b = \sum_{i} a_{i}\blacktriangleright b_{i}$. Because of the existence of local units in regular multiplier Hopf algebra,
 there is an $e\in A$ such that
 $ea_{i} = a_{i}$ for all $i = 1, 2, \ldots, n$. It follows easily that $e\blacktriangleright b = b$.
 So for any $b\in B$ there exists $e\in A$ such that $b = e\blacktriangleright b$.

 As an important consequence of the above result, we can use the Sweedler notation in the framework of dual pairs in the following sense.
 Take $a\in A$ and $b\in B$, and e.g. the element $b\blacktriangleright a = \<a_{(2)}, b\>a_{(1)}$.
 In the right hand side the element $a_{(2)}$ is covered by $b$ through the pairing because $b = e\blacktriangleright b$
 for some $e\in A$ and therefore $b\blacktriangleright a = \<a_{(2)}, b\> a_{(1)} = \<a_{(2)}, e\blacktriangleright b\> a_{(1)} = \<a_{(2)}e, b\> a_{(1)} \in A$.

 We also mention that for $a\in A$ and $b\in B$, $\<S(a), b\> = \<a, S(b)\>$, $\<a, 1\> = \varepsilon(a)$, $\<1, b\> = \varepsilon(b)$.
 For these formulas, one has to extend the pairing to $A\times M(B)$ and to $M(A)\times B$ (see Section 2 in \cite{D03a}).
 This can be done in a natural way using the fact that the four modules $A\blacktriangleright B$, $B\blacktriangleright A$,
 $A\blacktriangleleft B$ and $B\blacktriangleleft A$ are unital.
 Take e.g. $\<M(A), B\>$ as an example, we consider $\<x, b\>$, where $x\in M(A)$ and $b\in B$.
 For $b\in B$, there exists $e\in A$ such that $b = e\blacktriangleright b$.
 So $\<x, b\> = \<x, e\blacktriangleright b\> =\<xe, b\>$.
 \\

 A paring of two regular multiplier Hopf algebras is a natural setting for the construction of Drinfel'd double.
 It turns out that the conditions on the pairing $\<A, B\>$ are sufficient to construct the Drinfel'd double on $A \otimes B$.
 The main point of the essential ideals is that there is an invertible twist map $T: B \otimes A \rightarrow A \otimes B$,
 which defines an associative product on $A \otimes B$. For $a\in A$ and $b \in B$
 \begin{eqnarray}
 T(b \otimes a) = b_{(1)} \blacktriangleright a \blacktriangleleft S^{-1}(b_{(3)}) \otimes b_{(2)}. \label{2.1}
 \end{eqnarray}
 This map can be considered as a special case in \cite{YW11a}, and is bijective, and the inverse is given by
 \begin{eqnarray*}
 T^{-1}(a \otimes b) = b_{(2)} \otimes S^{-1}(b_{(1)}) \blacktriangleright a \blacktriangleleft b_{(3)}.
 \end{eqnarray*}

 As shown in \cite{DV04}, the structures of the Drinfel'd double is given as follows.

 Let $\mathscr{D} = A \bowtie B$ denote the algebra with tensor product $A \otimes B$ as the underlying space,
 and with the product given by the twist map $T$ as follows:
 \begin{eqnarray*}
 (a \bowtie b)(a' \bowtie b') = (m_{A} \otimes m_{B})(\iota \otimes T \otimes \iota) (a \otimes b \otimes a' \otimes b')
 \end{eqnarray*}
 with $a, a' \in A$ and $b, b' \in B$.
 If we write $T^{-1}(a' \otimes b') = \sum b'_{i} \otimes a'_{i}$, then we have
 \begin{eqnarray}
 (a \bowtie b)(a' \bowtie b') = \sum (a \otimes 1)T(b b'_{i} \otimes a'_{i}).
 \end{eqnarray}
 Similarly, if we write $T^{-1}(a \otimes b) = \sum b_{i} \otimes a_{i}$, then we have
 \begin{eqnarray}
 (a \bowtie b)(a' \bowtie b') = \sum T(b_{i} \otimes a_{i}a')(1 \otimes b'). \label{2.2}
 \end{eqnarray}
 There are algebra embeddings $A \rightarrow M(\mathscr{D}): \, a \rightarrow a \bowtie 1$ and
 $B \rightarrow M(\mathscr{D}): \, b \rightarrow 1 \bowtie b$. These embeddings can be extended to the multiplier algebras.

 The coproduct (or comultiplication), counit, and antipode are given as follows.
 \begin{eqnarray*}
 && \Delta (a \bowtie b) = \Delta_{A}^{cop}(a) \Delta_{B}(b), \, \,
  \varepsilon (a \bowtie b) = \varepsilon_{A} (a) \varepsilon_{B} (b), \, \,
  S(a \otimes b) = T(S_{B}(b) \otimes S^{-1}_{A}(a)).
 \end{eqnarray*}
 Without confusion, we always e.g. denote $S_{A}(a)$ just as $S(a)$.

 Drinfel'd double $\mathscr{D}$ can be considered as a special case of twisted double defined in \cite{YW11a}.

 \subsection{Complete modules}

 Let $A$ be a regular multiplier Hopf algebra. Suppose $X$ is a left $A$-module
 with the module structure map $\cdot : A\otimes X\longrightarrow  X$. We will always assume that the module is non-degenerate, this means
 that $x=0$ if $x\in X$ and $a\cdot x=0$ for all $a\in A$.
 If the module is unital (i.e., $A\cdot X=X$), then we can get an extension of the module structure to $M(A)$, this
 means that we can define $f\cdot x$, where $f\in M(A)$ and $x\in X$.
 In fact, since $x\in X=A\cdot X$, then $x= \sum_{i} a_{i}\cdot x_{i}$ and $f\cdot v=\sum_{i} (fa_{i})\cdot x_{i}$.
 In this setting, we can easily get $1_{M(A)}\cdot x=x$.

 Let $X$ be a left $A$-module. Denote by $Y$ the space of linear maps $\rho: A\rightarrow X$ satisfying
 $\rho(aa')=a\cdot \rho(a')$ for all $a, a'\in A$. Then $Y$ becomes a left $A$-module if we define $a\cdot \rho$ for
 $a\in A$ and $\rho \in Y$ by $(a\cdot \rho)(a')=\rho(a'a)=a'\cdot \rho(a)$.
 Define $\rho_{x} \in Y$ by $\rho_{x}(a)=a\cdot x$
 when $a\in A$. Then $X$ becomes a submodule of $Y$. Then we have $A\cdot Y\subseteq X$, and if $A\cdot X=X$, then
 $A\cdot Y=X$. Since $A^{2}=A$, $Y$ is still non-degenerate. If $A$ has a unit, then $Y=X$, in the other case,
 mostly $Y$ is strictly bigger than $X$. We can do the same as before for right modules as well.

 Let $X$ be a non-degenerate $A$-bimodule. Denotes by $Z$ the space of pair $(\lambda, \rho)$ of linear maps
 from $A$ to $X$ satisfying $a\cdot \lambda(a')=\rho(a)\cdot a'$ for all $a, a' \in A$. From the non-degeneracy, it follows
 that $\rho(aa')=a\cdot \rho(a')$ and $\lambda(aa')=\lambda(a)\cdot a'$ for all $a, a' \in A$. Also $\rho$ is completely
 determined by $\lambda$ and vice versa. We can consider $Z$ as the intersection of two extensions of $X$
 (as a left and a right modules). Then $Z$ becomes an $A$-bimodule, if we define $a \cdot z$ and $z\cdot a$
 for $a\in A$ and $z=(\lambda, \rho) \in Z$ by
 $a\cdot z = (a\lambda(\cdot), \rho(\cdot a))$ and $z\cdot a = (\lambda(a \cdot), \rho(\cdot) a)$. If we define
 $(\lambda_{x}, \rho_{x})$ for $x\in X$ by $\lambda_{x}(a)=x\cdot a$ and $\rho_{x}(a)=a\cdot x$,
 we get $X$ as a submodule of $Z$.
 We say that $Z$ is a completed module of $A$, and denote it as $M_{0}(X)$, see \cite{V08}.

 Let $V$ be a vector space and $X = A\otimes V$, we consider left and
 right action of $A$ on $X$ by for $a, a'\in A$ and $v\in V$,
 \begin{eqnarray*}
 a\cdot (a' \otimes v)= aa'\otimes v, \qquad
 (a' \otimes v)\cdot a= a'a\otimes v.
 \end{eqnarray*}
 The completed module we get here is denoted as $M_{0}(A\otimes V)$.

 \subsection{Yetter-Drinfel'd modules over a multiplier Hopf algebra}

 Recall of the definition of generalized (left-left) Yetter-Drinfel'd module over a multiplier Hopf algebra from \cite{YW11b}.

 Let $(A, \Delta, \varepsilon, S)$ be a regular multiplier Hopf algebra
 and $V$ a vector space. Then $V$ is called a Yetter-Drinfel'd module over $A$, if
 the following conditions hold:
 \begin{itemize}
 \item[(1)] $(V, \cdot)$ is an unital left $A$-module, i.e., $A\cdot V = V$;

 \item[(2)] $(V, \Gamma)$ is a left $A$-comodule, where $\Gamma: V\rightarrow M_{0}(A \otimes V)$ denotes the left coaction of $A$ on $V$;

 \item[(3)] $\Gamma$ and $\cdot$ satisfy the following condition:
 \begin{eqnarray}
 (a_{(1)}\cdot v)_{(-1)} a_{(2)}a' \otimes (a_{(1)}\cdot v)_{(0)} = a_{(1)}v_{(-1)}a' \otimes a_{(2)}\cdot v_{(0)} \label{2.3}
 \end{eqnarray}
 for all $a, a'\in A$ and $v\in V$.
 \end{itemize}

 By the definition of (left-left) Yetter-Drinfel'd modules, we can define Yetter-Drinfel'd module categories
 ${}_{A}^{A}\mathcal {YD}$. The objects in ${}_{A}^{A}\mathcal {YD}$ are left-left Yetter-Drinfel'd modules,
 and the morphisms are linear maps which interwine with the left action and the left coaction of $A$ on $M$,
 i.e., the morphisms between two objects are left $A$-linear and left $A$-colinear maps.
 Precisely, let $V, W \in {}_{A}^{A}\mathcal {YD}$ and $f: V\rightarrow W$ be a morphism, then
 \begin{eqnarray*}
 && f(a\cdot v) = a\cdot f(v), \\
 && (a' \otimes 1)\Gamma_{W}\circ f(v) = (a' \otimes 1) (\iota \otimes f)\Gamma_{V}(v),
 \end{eqnarray*}
 for all $a, a' \in A$ and $v \in V$, where $\Gamma_{W}$ ($\Gamma_{V}$) is the left coaction on $W$ ($V$).

 The other three kind of Yetter-Drinfel'd module categories are also defined in \cite{YW11b}.

\section{Heisenberg Double}
\def\theequation{\thesection.\arabic{equation}}
\setcounter{equation}{0}

 Let $\<A, B\>$ be a pairing of two regular multiplier Hopf algebras $A$ and $B$.
 The \emph{Heisenberg double} $\mathscr{H}$ is the smash product $A \# B$ with respect to the left regular action
 $b \blacktriangleright a = \<a_{(2)}, b\> a_{(1)}$ of $B$ on $A$. The production in $\mathscr{H}$ is given by
 \begin{eqnarray}
 && (a \# b)(a' \# b') = a(b_{(1)} \blacktriangleright a') \# b_{(2)}b'. \label{3.1}
 \end{eqnarray}
 for any $a, a' \in A$ and $b, b' \in B$.

 This definition can be found in \cite{DV04} and \cite{DrVZ99}.
 The product defined by (\ref{3.1}) is a twisted tensor product (see \cite{D03}) with a bijective twisted map
 $R(b \otimes a') = b_{(1)} \blacktriangleright a' \otimes b_{(2)}$, and the inverse is given by
 $R^{-1}(a \otimes b) = b_{(2)} \otimes S^{-1}(b_{(1)}) \blacktriangleright a$.
 By the Proposition 1.1 in \cite{D03}, this product is non-degenerate.
 \\

 Recall from \cite{S11} that for a finite dimensional Hopf algebra $B$ and its dual $B^{*}$,
 there is an action of Drinfel'd double $B^{*} \bowtie B$ on the Heisenberg double $B^{*} \# B$
 \begin{eqnarray*}
 && (f \bowtie b)\cdot (f' \# b') = f_{(3)} (b_{(1)} \blacktriangleright f') S^{-1}(f_{(2)}) \#
 (b_{(2)} b' S(b_{(3)})) \blacktriangleleft S^{-1}(f_{(1)}),
 \end{eqnarray*}
 where $f \bowtie b \in B^{*} \bowtie B$ and $f' \# b' \in B^{*} \# B$.
 In the following, we claim that this action also holds in the general multiplier Hopf algebra case.
 \\

 \textbf{Lemma 3.1}
 Let $\mathscr{D}$ be the Drinfel'd double and $\mathscr{H}$ the Heisenberg double for a multiplier Hopf algebra paring $\<A, B\>$,
 then for $a \bowtie b \in \mathscr{D}$ and $a' \# b' \in \mathscr{H}$,
 \begin{eqnarray}
 && (a \bowtie b)\cdot (a' \# b') = a_{(3)} (b_{(1)} \blacktriangleright a') S^{-1}(a_{(2)}) \#
 (b_{(2)} b' S(b_{(3)})) \blacktriangleleft S^{-1}(a_{(1)}) \label{3.2}
 \end{eqnarray}
 is well-difined.

 \emph{Proof}
 We need to check that (\ref{3.2}) makes sense in the framework of multiplier Hopf algebra paring.
 Indeed, for $a'\in A$ there is an $e \in B$ such that $e \blacktriangleright a' = a'$, then the right hand side
 \begin{eqnarray*}
 && a_{(3)} (b_{(1)} \blacktriangleright a') S^{-1}(a_{(2)}) \# (b_{(2)} b' S(b_{(3)})) \blacktriangleleft S^{-1}(a_{(1)}) \\
 &=& a_{(3)} (b_{(1)}e \blacktriangleright a') S^{-1}(a_{(2)}) \# (b_{(2)} b' S(b_{(3)})) \blacktriangleleft S^{-1}(a_{(1)}).
 \end{eqnarray*}
 $b_{(1)}, b_{(2)}$ are covered by $e$ and $b'$ respectively,
 $b_{(1)}e \otimes b_{(2)} b' \otimes b_{(3)} \in B\otimes B\otimes B$, so
 $b_{(1)} \blacktriangleright a' \otimes b_{(2)} b' S(b_{(3)}) \in A\otimes B$, it can be written as a finite sum of tensor products
 $\sum_{i} p_{i}\otimes q_{i}$ for some $p_{i} \in A$ and $q_{i} \in B$,
 then the right hand side is equal to
 $$\sum_{i} a_{(3)} p_{i} S^{-1}(a_{(2)}) \# q_{i} \blacktriangleleft S^{-1}(a_{(1)}).$$
 For $q_{i}$, there is an $f\in A$ so that $q_{i} = q_{i} \blacktriangleleft f$,
 so $a_{(1)}, a_{(3)}$ are covered by $S(f)$ and $p_{i}$ respectively, and the right hand side belongs to $A\otimes B$.
 $\hfill \Box$
 \\

 \textbf{Proposition \thesection.2}
 Let $\<A, B\>$ be a multiplier Hopf algebra paring, then by the action (\ref{3.2}), $\mathscr{H}$ is a unital $\mathscr{D}$-module.

 \emph{Proof}
 First we check that $\mathscr{H}$ is a $\mathscr{D}$-module, i.e., to show that
 $((c \bowtie d)(a \bowtie b)) \cdot (a' \# b') = (c \bowtie d) \cdot ((a \bowtie b) \cdot (a' \# b'))$.
 Indeed,
 \begin{eqnarray*}
 && (c \bowtie d) \cdot ((a \bowtie b)\cdot (a' \# b')) \\
 &=& (c \bowtie d) \cdot \Big( a_{(3)} (b_{(1)} \blacktriangleright a') S^{-1}(a_{(2)})
     \# (b_{(2)} b' S(b_{(3)})) \blacktriangleleft S^{-1}(a_{(1)}) \Big) \\
 &=& c_{(3)} \Big( d_{(1)} \blacktriangleright \big( a_{(3)} (b_{(1)} \blacktriangleright a') S^{-1}(a_{(2)}) \big) \Big) S^{-1}(c_{(2)}) \\
  && \# \Big( d_{(2)} \big( (b_{(2)} b' S(b_{(3)})) \blacktriangleleft S^{-1}(a_{(1)}) \big) S^{-1}(d_{(3)}) \Big)
     \blacktriangleleft S^{-1}(c_{(1)}) \\
 &=& c_{(3)} (d_{(1)} \blacktriangleright a_{(3)}) (d_{(2)}b_{(1)} \blacktriangleright a')
     (d_{(3)} \blacktriangleright S^{-1}(a_{(2)})) S^{-1}(c_{(2)}) \\
  && \# \Big( d_{(4)} \big( (b_{(2)} b' S(b_{(3)})) \blacktriangleleft S^{-1}(a_{(1)}) \big) S^{-1}(d_{(5)}) \Big)
     \blacktriangleleft S^{-1}(c_{(1)}) \\
 &=& c_{(3)} (d_{(1)} \blacktriangleright a_{(3)}) (d_{(2)}b_{(1)} \blacktriangleright a')
     S^{-1}(a_{(2)(2)}) S^{-1}(c_{(2)}) \\
  && \# \Big( (d_{(3)} \blacktriangleleft S^{-1}(a_{(2)(1)}))
     \big( (b_{(2)} b' S(b_{(3)})) \blacktriangleleft S^{-1}(a_{(1)}) \big) S^{-1}(d_{(4)}) \Big)
     \blacktriangleleft S^{-1}(c_{(1)}) \\
 &=& c_{(3)} (d_{(1)} \blacktriangleright a_{(4)}) (d_{(2)}b_{(1)} \blacktriangleright a')
     S^{-1}(a_{(3)}) S^{-1}(c_{(2)}) \\
  && \# \Big( (d_{(3)} \blacktriangleleft S^{-1}(a_{(2)}))
     \big( (b_{(2)} b' S(b_{(3)})) \blacktriangleleft S^{-1}(a_{(1)}) \big) S^{-1}(d_{(4)}) \Big)
     \blacktriangleleft S^{-1}(c_{(1)}),
 \end{eqnarray*}
 and
 \begin{eqnarray*}
 && \Big( (c \bowtie d)(a \bowtie b) \Big) \cdot (a' \# b') \\
 &=& \Big( c(d_{(1)} \blacktriangleright a \blacktriangleleft S^{-1}(d_{(3)})) \bowtie d_{(2)}b \Big) \cdot (a' \# b') \\
 &=& \Big(c \big( d_{(1)} \blacktriangleright a \blacktriangleleft S^{-1}(d_{(3)})\big) \Big)_{(3)}
     \Big((d_{(2)}b)_{(1)} \blacktriangleright a' \Big)
     S^{-1} \Big((c \big( d_{(1)} \blacktriangleright a \blacktriangleleft S^{-1}(d_{(3)})\big) )_{(2)} \Big) \\
  && \# \Big((d_{(2)}b)_{(2)} b' S \big( (d_{(2)}b)_{(3)} \big) \Big) \blacktriangleleft
     S^{-1} \Big((c \big( d_{(1)} \blacktriangleright a \blacktriangleleft S^{-1}(d_{(3)}) \big))_{(1)} \Big) \\
 &=& \Big(c_{(3)} \big(d_{(1)} \blacktriangleright a_{(3)} \big) \Big)
     \Big((d_{(2)}b)_{(1)} \blacktriangleright a' \Big)
     S^{-1} \Big(c_{(2)} a_{(2)} \Big) \\
  && \# \Big((d_{(2)}b)_{(2)} b' S \big( (d_{(2)}b)_{(3)} \big) \Big) \blacktriangleleft
     S^{-1} \Big( c_{(1)} (a_{(1)} \blacktriangleleft S^{-1}(d_{(3)})) \Big) \\
 &=& c_{(3)} \big(d_{(1)} \blacktriangleright a_{(3)} \big)
     \Big((d_{(2)(1)}b_{(1)}) \blacktriangleright a' \Big)
     S^{-1}(a_{(2)}) S^{-1}(c_{(2)}) \\
  && \# \Big(d_{(2)(2)}b_{(2)} b' S(b_{(3)})) S(d_{(2)(3)}) \Big) \blacktriangleleft
     \Big( S^{-1}(a_{(1)} \blacktriangleleft S^{-1}(d_{(3)})) S^{-1}(c_{(1)}) \Big) \\
 &=& c_{(3)} \big(d_{(1)} \blacktriangleright a_{(3)} \big)
     \Big((d_{(2)(1)}b_{(1)}) \blacktriangleright a' \Big)
     S^{-1}(a_{(2)}) S^{-1}(c_{(2)}) \\
  && \# \Big( \big(d_{(2)(2)} \blacktriangleleft S^{-1}(a_{(1)(2)(3)}) \big)
     \big( (b_{(2)} b' S(b_{(3)})) \blacktriangleleft S^{-1}(a_{(1)(2)(2)}) \big) \\
  && \big( S(d_{(2)(3)}) \blacktriangleleft S^{-1}(a_{(1)(2)(1)}) \big) \<a_{(1)(1)}, S^{-1}(d_{(3)})\>
     \Big) \blacktriangleleft S^{-1}(c_{(1)}) \\
 &=& c_{(3)} (d_{(1)} \blacktriangleright a_{(4)}) (d_{(2)}b_{(1)} \blacktriangleright a')
     S^{-1}(a_{(3)}) S^{-1}(c_{(2)}) \\
  && \# \Big( (d_{(3)} \blacktriangleleft S^{-1}(a_{(2)}))
     \big( (b_{(2)} b' S(b_{(3)})) \blacktriangleleft S^{-1}(a_{(1)}) \big) S^{-1}(d_{(4)}) \Big)
     \blacktriangleleft S^{-1}(c_{(1)}).
 \end{eqnarray*}

 Then, we will show that the module action is unital. We denote the adjoint actions of $A$ and $B$ on themselves by
 $a\rightharpoonup a' = a_{(2)} a' S^{-1}(a_{(1)})$, $b \rightharpoondown b' = b_{(1)} b' S(b_{(2)})$,
 it is easy to show that these two actions are unital. Define $F, G: A\otimes B\rightarrow A\otimes B$
 by $F(a \otimes b)=a_{(2)} \otimes b\blacktriangleleft S^{-1}(a_{(1)})$
 and $G(a \otimes b)=b_{(1)}\blacktriangleright a \otimes b_{(2)}$, $F$ and $G$ are bijective. So
 \begin{eqnarray*}
 && (a \bowtie b)\cdot (a' \# b')
 = (\rightharpoonup\otimes\iota)F_{13}(\iota\otimes\iota\otimes\rightharpoondown)(\iota\otimes G \otimes\iota)(a\otimes a'\otimes b\otimes b'),
 \end{eqnarray*}
 this can conclude that the action of $\mathscr{D}$ on $\mathscr{H}$ is unital.
 $\hfill \Box$
 \\

 \textbf{Remark}
 Since $\mathscr{H}$ is a unital $\mathscr{D}$-module, we can get an extension of the module structure to $M(\mathscr{D})$,
 this means we can define $f \cdot (a\# b)$, where $f\in M(\mathscr{D})$ and $a\# b \in \mathscr{H}$.
 Indeed, since $a\# b \in \mathscr{H} = \mathscr{D} \cdot \mathscr{H}$, then $a\# b = \sum_{i, j}(a'_{j}\bowtie b'_{j})\cdot(a_{i}\# b_{i})$,
 $f \cdot (a\# b) = \sum_{i, j} \big(f (a'_{j}\bowtie b'_{j})\big) \cdot(a_{i}\# b_{i})$.
 \\

 \textbf{Proposition \thesection.3}
 $\mathscr{H}$ is a $\mathscr{D}$-module algebra.

 \emph{Proof}
 It is sufficient to show that
 \begin{eqnarray*}
 (a \bowtie b)\cdot \Big((a' \# b')(c' \# d') \Big) =
 \Big( (a_{(2)} \bowtie b_{(1)}) \cdot(a' \# b') \Big)\Big( (a_{(1)} \bowtie b_{(2)}) \cdot (c' \# d') \Big).
 \end{eqnarray*}
 In fact,
 \begin{eqnarray*}
 && (a \bowtie b)\cdot \Big((a' \# b')(c' \# d') \Big) \\
 &=& (a \bowtie b)\cdot \Big(a'(b'_{(1)} \blacktriangleright c') \# b'_{(2)} d' \Big) \\
 &=& a_{(3)} \Big(b_{(1)} \blacktriangleright (a'(b'_{(1)} \blacktriangleright c')) \Big) S^{-1}(a_{(2)})
     \# \Big(b_{(2)} (b'_{(2)} d') S(b_{(3)}) \Big) \blacktriangleleft S^{-1}(a_{(1)}) \\
 &=& a_{(3)} (b_{(1)} \blacktriangleright a')\big((b_{(2)}b'_{(1)}) \blacktriangleright c' \big)  S^{-1}(a_{(2)})
     \# \Big(b_{(3)} b'_{(2)} d' S(b_{(4)}) \Big) \blacktriangleleft S^{-1}(a_{(1)}),
 \end{eqnarray*}
 and
 \begin{eqnarray*}
 && \Big( (a_{(2)} \bowtie b_{(1)}) \cdot(a' \# b') \Big)\Big( (a_{(1)} \bowtie b_{(2)}) \cdot (c' \# d') \Big) \\
 &=& \Big(a_{(6)} (b_{(1)}\blacktriangleright a') S^{-1}(a_{(5)}) \# (b_{(2)} b' S(b_{(3)})) \blacktriangleleft S^{-1}(a_{(4)}) \Big) \\
  && \Big(a_{(3)} (b_{(4)}\blacktriangleright c') S^{-1}(a_{(2)}) \# (b_{(5)} d' S(b_{(6)})) \blacktriangleleft S^{-1}(a_{(1)}) \Big) \\
 &=& \Big(a_{(4)} (b_{(1)}\blacktriangleright a') \big((b_{(2)} b' S(b_{(3)}))_{(1)} \blacktriangleright S^{-1}(a_{(3)})\big)
          \# (b_{(2)} b' S(b_{(3)}))_{(2)} \Big) \\
  && \Big(a_{(2)} (b_{(4)}\blacktriangleright c') \big((b_{(5)} d' S(b_{(6)}))_{(1)} \blacktriangleright S^{-1}(a_{(1)})\big)
          \# (b_{(5)} d' S(b_{(6)}))_{(2)} \Big) \\
 &=& a_{(4)} (b_{(1)}\blacktriangleright a') \Big((b_{(2)} b' S(b_{(3)}))_{(1)} \blacktriangleright S^{-1}(a_{(3)})\Big)
         \Big((b_{(2)} b' S(b_{(3)}))_{(2)} \blacktriangleright \big(a_{(2)} (b_{(4)}\blacktriangleright c') \big)\Big) \\
  && \Big((b_{(2)} b' S(b_{(3)}))_{(3)} (b_{(5)} d' S(b_{(6)}))_{(1)} \blacktriangleright S^{-1}(a_{(1)})\Big)
     \# (b_{(2)} b' S(b_{(3)}))_{(4)} (b_{(5)} d' S(b_{(6)}))_{(2)} \\
 &=& a_{(4)} (b_{(1)}\blacktriangleright a') \Big((b_{(2)} b' S(b_{(3)}))_{(1)} \blacktriangleright
     \big(S^{-1}(a_{(3)})a_{(2)} (b_{(4)}\blacktriangleright c') \big)\Big) \\
  && \Big((b_{(2)} b' S(b_{(3)}))_{(2)} (b_{(5)} d' S(b_{(6)}))_{(1)} \blacktriangleright S^{-1}(a_{(1)})\Big)
     \# (b_{(2)} b' S(b_{(3)}))_{(3)} (b_{(5)} d' S(b_{(6)}))_{(2)} \\
 &=& a_{(3)} (b_{(1)} \blacktriangleright a')\big((b_{(2)}b'_{(1)}) \blacktriangleright c' \big)  S^{-1}(a_{(2)})
     \# \Big(b_{(3)} b'_{(2)} d' S(b_{(4)}) \Big) \blacktriangleleft S^{-1}(a_{(1)}).
 \end{eqnarray*}
 This completes the proof.
 $\hfill \Box$
 \\

 Now, we give $\mathscr{H}$ a coaction as follows: $\Gamma: \mathscr{H} \longrightarrow M_{0}(\mathscr{D} \otimes \mathscr{H})$,
 \begin{eqnarray}
 ((a'\bowtie b')\otimes 1)\Gamma(a \# b) = (a'\bowtie b')(a_{(2)}\bowtie b_{(1)}) \otimes a_{(1)}\# b_{(2)}, \label{3.3}
 \end{eqnarray}
 where $M_{0}(\mathscr{D} \otimes \mathscr{H})$ is the complete module (see \cite{DVW11,YW11b}).

 Obviously, this coaction is well-defined.
 \\

 \textbf{Proposition \thesection.4}
 The coaction $\Gamma$ defined above makes $\mathscr{H}$ a left $\mathscr{D}$-comodule algebra.

 \emph{Proof}
 We firstly can check that $\Gamma$ is well-defined by the formula (\ref{2.2}),
 and satisfies $(\Delta\otimes\iota)\Gamma=(\iota\otimes\Gamma)\Gamma$.
 Then, we will show that $\Gamma$ defined above is injective.
 Indeed, if $\Gamma(a \# b)=0$, applying $\varepsilon\otimes\iota$ on this equation, we can get that $a \# b=0$,
 so it is injective.

 Finally, we need to show $\Gamma$ satisfy $\Gamma((a' \# b')(c' \# d')) = \Gamma(a' \# b')\Gamma(c' \# d')$.
 Indeed,
 \begin{eqnarray*}
 && \big((a'\bowtie b')\otimes 1 \big)\Gamma\big((a \# b)(c \# d)\big)  \\
 &=& \big((a'\bowtie b')\otimes 1 \big)\Gamma\big(a(b_{(1)} \blacktriangleright c) \# b_{(2)}d \big)  \\
 &=& \big((a'\bowtie b')\otimes 1 \big)\Big(a(b_{(1)} \blacktriangleright c)\Big)_{(2)} \bowtie (b_{(2)}d)_{(1)}
     \otimes \Big(a(b_{(1)} \blacktriangleright c)\Big)_{(1)} \# (b_{(2)}d)_{(2)} \\
 &=& (a'\bowtie b')\Big(a_{(2)}(b_{(1)} \blacktriangleright c_{(2)}) \bowtie b_{(2)}d_{(1)}\Big)
     \otimes a_{(1)} c_{(1)} \# b_{(3)}d_{(2)},
 \end{eqnarray*}
 and
 \begin{eqnarray*}
 && \big((a'\bowtie b')\otimes 1 \big)\big(\Gamma(a \# b)\Gamma(c \# d)\big)  \\
 &=& (a'\bowtie b')(a_{(2)}\bowtie b_{(1)})(c_{(2)}\bowtie d_{(1)}) \otimes (a_{(1)} \# b_{(2)})(c_{(1)} \# d_{(2)}) \\
 &=& (a'\bowtie b')
     \Big(a_{(2)}\big( b_{(1)(1)}\blacktriangleright c_{(2)} \blacktriangleleft S^{-1}(b_{(1)(3)}) \big) \bowtie b_{(1)(2)}d_{(1)}\Big) \\
  && \otimes a_{(1)}(b_{(2)(1)} \blacktriangleright c_{(1)}) \# b_{(2)(2)} d_{(2)} \\
 &=& (a'\bowtie b')
     \Big(a_{(2)}\big( b_{(1)(1)}\blacktriangleright (c_{(2)}\blacktriangleleft b_{(2)(1)}) \blacktriangleleft S^{-1}(b_{(1)(3)}) \big) \bowtie b_{(1)(2)}d_{(1)}\Big) \\
  && \otimes a_{(1)} c_{(1)} \# b_{(2)(2)} d_{(2)} \\
 &=& (a'\bowtie b')\Big(a_{(2)}\big(b_{(1)}\blacktriangleright c_{(2)}\big)\bowtie b_{(2)}d_{(1)}\Big) \otimes a_{(1)}c_{(1)} \# b_{(3)}d_{(2)}.
 \end{eqnarray*}
 This completes the proof.
 $\hfill \Box$
 \\

 By a Yetter-Drinfel'd module algebra, we mean a module and comodule algebra that is also a Yetter-Drinfel'd module
 (see Section 2.2 or paper \cite{YW11b}).
 Then we can get the first main result of this paper.
 \\

 \textbf{Theorem \thesection.5}
 For the Drinfel'd double $\mathscr{D}$ and Heisenberg double $\mathscr{H}$ based on a pairing of regular multiplier Hopf algebras $A$ and $B$,
 $\mathscr{H}$ endowed with action (\ref{3.2}) and coaction (\ref{3.3}) is a (left-left) $\mathscr{D}$-Yetter-Drinfel'd module algebra.

 \emph{Proof}
 From Proposition \thesection.4 and \thesection.5, we know that $\mathscr{H}$ is a $\mathscr{D}$-module and comodule algebra.
 the left thing we need to do is checking the compatible condition (\ref{2.3}) of Yetter-Drinfel'd module
 ${}_{\mathscr{D}}^{\mathscr{D}}\mathcal {YD}$, i.e.,
 \begin{eqnarray*}
 && \Big((a \bowtie b)_{(1)} \cdot (c \# d)\Big)_{(-1)} (a \bowtie b)_{(2)} (a'\bowtie b')
    \otimes \Big((a \bowtie b)_{(1)} \cdot (c \# d)\Big)_{(0)} \\
 &=& (a \bowtie b)_{(1)} (c \# d)_{(-1)} (a'\bowtie b') \otimes (a \bowtie b)_{(2)} \cdot (c \# d)_{(0)}.
 \end{eqnarray*}
 Indeed,
 \begin{eqnarray*}
 && (a \bowtie b)_{(1)} (c \# d)_{(-1)} (a'\bowtie b') \otimes (a \bowtie b)_{(2)} \cdot (c \# d)_{(0)} \\
 &=& (a_{(2)} \bowtie b_{(1)}) (c_{(2)} \bowtie d_{(1)}) (a'\bowtie b') \otimes (a_{(1)} \bowtie b_{(2)}) \cdot (c_{(1)} \# d_{(2)}) \\
 &=& (a_{(4)} \bowtie b_{(1)}) \Big(c_{(2)}(d_{(1)} \blacktriangleright a' \blacktriangleleft S^{-1}(d_{(3)})) \bowtie d_{(2)}b'\Big) \\
  && \otimes a_{(3)} (b_{(2)}\blacktriangleright c_{(1)})S^{-1}(a_{(2)}) \# \big(b_{(3)}d_{(4)}S(b_{(4)})\big)\blacktriangleleft S^{-1}(a_{(1)}) \\
 &=& a_{(4)} \Big( b_{(1)} \blacktriangleright \big(c_{(2)}(d_{(1)} \blacktriangleright a' \blacktriangleleft S^{-1}(d_{(3)})) \big)
     \blacktriangleleft S^{-1}(b_{(3)}) \Big)  \bowtie b_{(2)}d_{(2)}b' \\
  && \otimes a_{(3)} (b_{(4)}\blacktriangleright c_{(1)})S^{-1}(a_{(2)}) \# \big(b_{(5)}d_{(4)}S(b_{(6)})\big)\blacktriangleleft S^{-1}(a_{(1)}) \\
 &=& a_{(4)} \big( b_{(1)(1)} \blacktriangleright c_{(2)} \blacktriangleleft S^{-1}(b_{(3)(2)}) \big)
     \big( b_{(1)(2)}d_{(1)} \blacktriangleright  a' \blacktriangleleft S^{-1}(d_{(3)}) S^{-1}(b_{(3)(1)}) \big) \\
  && \bowtie b_{(2)}d_{(2)}b'
     \otimes a_{(3)} (b_{(4)}\blacktriangleright c_{(1)})S^{-1}(a_{(2)})
     \# \big(b_{(5)}d_{(4)}S(b_{(6)})\big)\blacktriangleleft S^{-1}(a_{(1)}) \\
 &=& a_{(4)} \big( b_{(1)} \blacktriangleright c_{(2)} \big)
     \big( b_{(2)}d_{(1)} \blacktriangleright  a' \blacktriangleleft S^{-1}(b_{(4)} d_{(3)}) \big)
     \bowtie b_{(3)}d_{(2)}b' \\
  && \otimes a_{(3)} c_{(1)} S^{-1}(a_{(2)})
     \# \big(b_{(5)}d_{(2)}S(b_{(6)})\big) \blacktriangleleft S^{-1}(a_{(1)}) \\
 &=& \Big( a_{(4)} \big( b_{(1)} \blacktriangleright c_{(2)} \big) \bowtie b_{(2)}d_{(1)} \Big) (a'\bowtie b')  \\
  && \otimes a_{(3)} c_{(1)} S^{-1}(a_{(2)}) \# \big(b_{(3)}d_{(2)}S(b_{(4)})\big) \blacktriangleleft S^{-1}(a_{(1)}),
 \end{eqnarray*}
 and
 \begin{eqnarray*}
 && \Big((a \bowtie b)_{(1)} \cdot (c \# d)\Big)_{(-1)} (a \bowtie b)_{(2)} (a'\bowtie b')
     \otimes \Big((a \bowtie b)_{(1)} \cdot (c \# d)\Big)_{(0)} \\
 &=& \Big((a_{(2)} \bowtie b_{(1)}) \cdot (c \# d)\Big)_{(-1)} (a_{(1)} \bowtie b_{(2)}) (a'\bowtie b')
     \otimes \Big((a_{(2)} \bowtie b_{(1)}) \cdot (c \# d)\Big)_{(0)} \\
 &=& \Big(a_{(4)} (b_{(1)} \blacktriangleright c) S^{-1}(a_{(3)}) \# \big(b_{(2)} d S(b_{(3)}) \big) \blacktriangleleft S^{-1}(a_{(2)}) \Big)_{(-1)}
     (a_{(1)} \bowtie b_{(4)}) (a'\bowtie b') \\
  && \otimes \Big(a_{(4)} (b_{(1)} \blacktriangleright c) S^{-1}(a_{(3)}) \#
     \big(b_{(2)} d S(b_{(3)}) \big) \blacktriangleleft S^{-1}(a_{(2)}) \Big)_{(0)} \\
 &=& \Big(a_{(6)} (b_{(1)} \blacktriangleright c_{(2)}) S^{-1}(a_{(3)}) \bowtie \big(b_{(2)} d S(b_{(3)}) \big)_{(1)}
     \blacktriangleleft S^{-1}(a_{(2)}) \Big) (a_{(1)} \bowtie b_{(4)}) (a'\bowtie b') \\
  && \otimes \Big(a_{(5)}  c_{(1)} S^{-1}(a_{(4)}) \# \big(b_{(2)} d S(b_{(3)}) \big)_{(2)} \Big) \\
 &=& \Big(a_{(5)} (b_{(1)} \blacktriangleright c_{(2)}) \big((b_{(2)} d S(b_{(3)}))_{(1)} \blacktriangleright S^{-1}(a_{(2)})\big)
     \bowtie \big(b_{(2)} d S(b_{(3)}) \big)_{(2)} \Big) \\
  && (a_{(1)} \bowtie b_{(4)}) (a'\bowtie b')
     \otimes \Big(a_{(4)}  c_{(1)} S^{-1}(a_{(3)}) \# \big(b_{(2)} d S(b_{(3)}) \big)_{(3)} \Big) \\
 &=& \Big( a_{(5)} (b_{(1)} \blacktriangleright c_{(2)}) \big((b_{(2)} d S(b_{(3)}))_{(1)} \blacktriangleright S^{-1}(a_{(2)})\big) \\
  && \big( (b_{(2)} d S(b_{(3)}))_{(2)} \blacktriangleright a_{(1)} \blacktriangleleft S^{-1}((b_{(2)} d S(b_{(3)}))_{(4)}) \big)
     \bowtie (b_{(2)} d S(b_{(3)}))_{(3)} b_{(4)} \Big)(a'\bowtie b') \\
  && \otimes \Big(a_{(4)}  c_{(1)} S^{-1}(a_{(3)}) \# \big(b_{(2)} d S(b_{(3)}) \big)_{(5)} \Big) \\
 &=& \Big( a_{(6)} (b_{(1)} \blacktriangleright c_{(2)}) \big((b_{(2)} d S(b_{(3)}))_{(1)} \blacktriangleright S^{-1}(a_{(3)})\big) \\
  && \big( a_{(1)} \blacktriangleleft S^{-1}((b_{(2)} d S(b_{(3)}))_{(3)}) \big)
     \bowtie \big( (b_{(2)} d S(b_{(3)}))_{(2)} \blacktriangleleft a_{(2)}\big) b_{(4)} \Big)(a'\bowtie b') \\
  && \otimes \Big(a_{(5)}  c_{(1)} S^{-1}(a_{(4)}) \# \big(b_{(2)} d S(b_{(3)}) \big)_{(5)} \Big) \\
 &=& \Big( a_{(5)} (b_{(1)} \blacktriangleright c_{(2)})  S^{-1}(a_{(2)})
     \big( a_{(1)} \blacktriangleleft S^{-1}((b_{(2)} d S(b_{(3)}))_{(2)}) \big) \\
  && \bowtie (b_{(2)} d S(b_{(3)}))_{(1)} b_{(4)} \Big)(a'\bowtie b')
     \otimes \Big(a_{(4)}  c_{(1)} S^{-1}(a_{(3)}) \# \big(b_{(2)} d S(b_{(3)}) \big)_{(3)} \Big) \\
 &=& \Big( a_{(4)} \big( b_{(1)} \blacktriangleright c_{(2)} \big) \bowtie b_{(2)}d_{(1)} \Big) (a'\bowtie b')  \\
  && \otimes a_{(3)} c_{(1)} S^{-1}(a_{(2)}) \# \big(b_{(3)}d_{(2)}S(b_{(4)})\big) \blacktriangleleft S^{-1}(a_{(1)}).
 \end{eqnarray*}
 This completes the proof.
 $\hfill \Box$
 \\

 \textbf{Example \thesection.6}
 Let $G$ be a group with unit $e$, $B=K[G]$ the group Hopf algebra and $A=K(G)$ the well-known multiplier
 Hopf algebra on $G$. Then the Heisenberg double $\mathscr{D} = A\# B$ with the product
 $$(\delta_{p} \# q)(\delta_{p'} \# q') = \delta_{p} \delta_{p'q^{-1}} \# qq',
 $$
 and the Drinfel'd double $\mathscr{D} = A\bowtie B$ with the structures
 as follows. The product in $\mathscr{D}$ is given by
 $$
 (\delta_{g}\bowtie h)(\delta_{p}\bowtie q) = \delta_{g}\delta_{h ph^{-1}} \bowtie hq
 $$
 for all $\delta_{g}, \delta_{p} \in K(G), h, q\in K[G]$, and the multiplier Hopf structure is given by
 \begin{eqnarray*}
 && \Delta(\delta_{g}\otimes h) = \sum_{p\in G} (\delta_{p^{-1}g}\otimes h)\otimes (\delta_{p}\otimes h), \\
 && \varepsilon(\delta_{g} \otimes h) = \delta_{g, e},  \\
 && S(\delta_{g}\otimes h) = \delta_{h^{-1} g^{-1}h}\otimes h^{-1}.
 \end{eqnarray*}
 Then the action $(\delta_{p} \bowtie q) \cdot (\delta_{p'} \# q')=\delta_{p'q'q^{-1}p} \delta_{p'q^{-1}} \# qq'q^{-1}$
 and the coaction $\Gamma(\delta_{p} \# q) = \sum_{s\in G} \delta_{s^{-1}p} \bowtie q \otimes \delta_{s} \#q$
 make $K(G)\# K[G]$ a $K(G)\bowtie K[G]$-Yetter-Drinfel'd module algebra.
 \\

 Let $A$ be a regular multiplier Hopf algebra, a left $A$-module and left $A$-comodule algebra $X$ is said to be
 \emph{braided commutative} (or $A$-commutative in \cite{YZ13}), if for any $x, y \in X$
 \begin{eqnarray}
 yx=(y_{(-1)} \cdot x)y_{(0)}.
 \end{eqnarray}
 For any two (left-left) Yetter-Drinfel'd $A$-module algebras $X$ and $Y$, their \emph{braided product}
 (shown in the proof of Theorem 4.1 in \cite{D07})
 $X\propto Y$ is defined as follows
 \begin{eqnarray}
 (x\propto y)(x'\propto y') = x(y_{(-1)} \cdot x') \propto y_{(0)}y'.
 \end{eqnarray}
 for $x, x' \in X$ and $y, y' \in Y$.
 \\

 \textbf{Proposition \thesection.7}
 $X\propto Y$ is a Yetter-Drinfel'd $A$-module algebra.

 \emph{Proof} Let $t_{Y, X}(y\otimes x') = y_{(-1)} \cdot x' \otimes y_{(0)}$,
 then equation $(3.5)$ defines a twisted tensor product algebra in \cite{D03}.
 By Proposition 2.3 in \cite{D07},
 we can easily get $X\propto Y$ is an $A$-module and $A$-comodule algebra satisfying the compatibility condition
 of Yetter-Drinfel'd module, i.e., a Yetter-Drinfel'd $A$-module algebra.

 In detail, firstly we check that $X\propto Y$ is an $A$-module algebra, i.e.,
 $a\cdot ((x\propto y)(x'\propto y')) = (a_{(1)}\cdot (x\propto y))(a_{(2)}\cdot (x'\propto y'))$.
 Indeed, the $A$-module action on $x\propto y$ is given by $a\cdot (x\propto y) = (a_{(1)}\cdot x)\propto (a_{(2)}\cdot y)$, and
 \begin{eqnarray*}
 && (a_{(1)}\cdot (x\propto y))(a_{(2)}\cdot (x'\propto y')) \\
 &=& (a_{(1)}\cdot x\propto a_{(2)}\cdot y)(a_{(3)}\cdot x'\propto a_{(4)}\cdot y') \\
 &\stackrel{(3.5)}{=}& (a_{(1)}\cdot x)((a_{(2)}\cdot y)_{(-1)}a_{(3)}\cdot x') \propto (a_{(2)}\cdot y)_{(0)}(a_{(4)}\cdot y') \\
 &\stackrel{(2.4)}{=}& (a_{(1)}\cdot x)(a_{(2)} y_{(-1)} \cdot x') \propto (a_{(3)}\cdot y_{(0)})(a_{(4)}\cdot y') \\
 &=& a_{(1)}\cdot (x(y_{(-1)} \cdot x')) \propto a_{(1)} \cdot (y_{(0)}y') \\
 &=& a\cdot (x(y_{(-1)} \cdot x') \propto y_{(0)}y') \\
 &\stackrel{(3.5)}{=}& a\cdot ((x\propto y)(x'\propto y')).
 \end{eqnarray*}

 Secondly, $X\propto Y$ is an $A$-comodule algebra. We need to show
 $\Gamma((x\propto y)(x'\propto y')) = \Gamma(x\propto y) \Gamma(x'\propto y')$. For any $a\in A$,
 $\Gamma(x\propto y)(a\otimes 1) = x_{(-1)}y_{(-1)}a \otimes x_{(0)}\propto y_{(0)}$, then
 \begin{eqnarray*}
 && \Gamma((x\propto y)(x'\propto y'))(a\otimes 1) \\
 &=& \Gamma(x(y_{(-1)} \cdot x') \propto y_{(0)}y')(a\otimes 1) \\
 &=& (x(y_{(-1)} \cdot x'))_{(-1)} (y_{(0)}y')_{(-1)} a\otimes (x(y_{(-1)} \cdot x'))_{(0)} \propto (y_{(0)}y')_{(0)} \\
 &=& x_{(-1)} (y_{(-1)} \cdot x')_{(-1)} y_{(0)(-1)}y'_{(-1)} a\otimes x_{(0)}(y_{(-1)} \cdot x')_{(0)} \propto y_{(0)(0)}y'_{(0)} \\
 &\stackrel{(2.4)}{=}& x_{(-1)} y_{(-2)} x'_{(-1)} y'_{(-1)} a\otimes x_{(0)}(y_{(-1)} \cdot x'_{(0)}) \propto y_{(0)}y'_{(0)} \\
 &\stackrel{(3.5)}{=}& x_{(-1)} y_{(-1)} x'_{(-1)} y'_{(-1)} a\otimes (x_{(0)}\propto y_{(0)})(x'_{(0)} \propto y'_{(0)}) \\
 &=& \Gamma(x\propto y) \Gamma(x'\propto y')(a\otimes 1).
 \end{eqnarray*}

 Finally, we need to check the Yetter-Drinfel¡¯d compatibility condition.
 \begin{eqnarray*}
 && (a_{(1)}\cdot (x\propto y))_{(-1)} a_{(2)} a' \otimes (a_{(1)}\cdot (x\propto y))_{(0)} \\
 &=& ((a_{(1)}\cdot x)\propto (a_{(2)}\cdot y))_{(-1)} a_{(3)} a' \otimes ((a_{(1)}\cdot x)\propto (a_{(2)}\cdot y))_{(0)} \\
 &=& (a_{(1)}\cdot x)_{(-1)} (a_{(2)}\cdot y)_{(-1)} a_{(3)} a' \otimes (a_{(1)}\cdot x)_{(0)} \propto (a_{(2)}\cdot y)_{(0)} \\
 &\stackrel{(2.4)}{=}& (a_{(1)}\cdot x)_{(-1)} a_{(2)} y_{(-1)} a' \otimes (a_{(1)}\cdot x)_{(0)} \propto (a_{(3)}\cdot y_{(0)}) \\
 &\stackrel{(2.4)}{=}& a_{(1)} x_{(-1)}  y_{(-1)} a' \otimes (a_{(2)}\cdot x_{(0)}) \propto (a_{(3)}\cdot y_{(0)}) \\
 &=& a_{(1)} (x_{(-1)}  y_{(-1)}) a' \otimes a_{(2)}\cdot (x_{(0)} \propto  y_{(0)}) \\
 &=& a_{(1)} (x\propto y)_{(-1)} a' \otimes a_{(2)}\cdot (x\propto y)_{(0)}.
 \end{eqnarray*}
 This completes the proof. $\hfill \Box$
 \\

 Now, we can get another main result of this paper.

 \textbf{Theorem \thesection.8}
 $\mathscr{H}$ is a braided $\mathscr{D}$-commutative algebra.
 And $\mathscr{H}$ is the braided product $A \propto B$,
 where $A$ and $B$ are braided commutative Yetter-Drinfel'd $\mathscr{D}$-module algebras by restriction, i.e.,
 the action is given by
 \begin{eqnarray*}
 && (a \bowtie b)\cdot a' = a_{(2)} (b \blacktriangleright a') S^{-1}(a_{(1)}), \\
 && (a \bowtie b)\cdot b' = (b_{(1)} b' S(b_{(2)})) \blacktriangleleft S^{-1}(a),
 \end{eqnarray*}
 and coaction $\rho$: $A \rightarrow M_{0}(\mathscr{D} \otimes A)$ and $B \rightarrow M_{0}(\mathscr{D} \otimes B)$
 is given by $\rho(a')=\Delta^{cop}_{13}(a')$, $\rho(b')=\Delta_{23}(b')$,
 for $a \bowtie b \in \mathscr{D}$, $a' \in A$ and $b' \in B$.

 \emph{Proof}
 We need to show that $(a' \# b')(a \# b) = \big((a' \# b')_{(-1)} \cdot (a \# b) \big)(a' \# b')_{(0)}$.
 Indeed,
 \begin{eqnarray*}
 && \big((a' \# b')_{(-1)} \cdot (a \# b) \big)(a' \# b')_{(0)} \\
 &=& \big((a'_{(2)} \bowtie b'_{(1)}) \cdot (a \# b) \big)(a'_{(1)} \# b'_{(2)}) \\
 &=& \Big( a'_{(4)} (b'_{(1)} \blacktriangleright a) S^{-1}(a'_{(3)}) \# (b'_{(2)} b S(b'_{(3)})) \blacktriangleleft S^{-1}(a'_{(2)}) \Big)
     (a'_{(1)} \# b'_{(4)}) \\
 &=& \Big( a'_{(3)} (b'_{(1)} \blacktriangleright a) \big((b'_{(2)} b S(b'_{(3)}))_{(1)} \blacktriangleright S^{-1}(a'_{(2)}) \big)
     \# (b'_{(2)} b S(b'_{(3)}))_{(2)} \Big)  (a'_{(1)} \# b'_{(4)}) \\
 &=& a'_{(3)} (b'_{(1)} \blacktriangleright a) \big((b'_{(2)} b S(b'_{(3)}))_{(1)} \blacktriangleright S^{-1}(a'_{(2)}) \big)
     \Big( (b'_{(2)} b S(b'_{(3)}))_{(2)} \blacktriangleright  a'_{(1)} \Big) \\
  && \# (b'_{(2)} b S(b'_{(3)}))_{(3)} b'_{(4)} \\
 &=& a'_{(3)} (b'_{(1)} \blacktriangleright a) \Big((b'_{(2)} b S(b'_{(3)}))_{(1)} \blacktriangleright (S^{-1}(a'_{(2)})a'_{(1)}) \Big)
     \# (b'_{(2)} b S(b'_{(3)}))_{(3)} b'_{(2)} \\
 &=& a' (b'_{(1)} \blacktriangleright a) \# b'_{(2)} b \\
 &=& (a' \# b') (a\# b).
 \end{eqnarray*}
 This shows that $\mathscr{H}$ is a braided $\mathscr{D}$-commutative algebra.

 The second part is obvious, since $\mathscr{H}$ is a unital $\mathscr{D}$-module,
 we can get an extension of the module structure to $M(\mathscr{D})$, and
 \begin{eqnarray*}
 (a\propto b)(a'\propto b')
 &=& a(b_{(-1)} \cdot a') \propto b_{(0)}b' \\
 &=& a((1\bowtie b_{(1)}) \cdot a') \propto b_{(2)}b' \\
 &=& a(b_{(1)} \blacktriangleright a') \propto b_{(2)} b'.
 \end{eqnarray*}
 This completes the proof.
 $\hfill \Box$
 \\

 In the end of this section, we consider $\<\widehat{A}, A\>$, where $A$ is a regular multiplier Hopf algebra with a left integral $\varphi$,
 and $\widehat{A}= \varphi(\cdot A)$ be the duality introduced in \cite{V98}.

 In this situation, the conditions of a pairing on $\<\widehat{A}, A\>$ naturally hold.
 Indeed, it is easy to check that $\<\widehat{A}, A\>$ is a pre-pairing,
 we only need to check one of six equivalent conditions in Proposition 2.7 of \cite{DrV01}.
 Because for $a, b\in A$, $b\blacktriangleright \varphi(\cdot a) = \varphi(\cdot ba)$ and $A^{2}=A$,
 we get $A\blacktriangleright \widehat{A} = \widehat{A}$, i.e., condition (2) in Proposition 2.7 of \cite{DrV01} holds.
 So $\<\widehat{A}, A\>$ is a pairing, furthermore a special case of $\<A, B\>$ introduced before.
 \\

 \textbf{Corollary \thesection.9}
 Let $A$ be a regular multiplier Hopf algebra with a left integral $\varphi$, and $\widehat{A}$ be the dual regular multiplier Hopf algebra.
 Then Heisenberg double $\mathscr{H} = \widehat{A}\# A$ is a (left-left) $\mathscr{D} = \widehat{A}\bowtie A$-Yetter-Drinfel'd module algebra,
 and moreover $\mathscr{H}$ is a braided $\mathscr{D}$-commutative algebra.
 \\

 \textbf{Example \thesection.10}
 Take the notations as Example 3.6, $K(G)\# K[G]$ is a braided $K(G)\bowtie K[G]$-commutative algebra.
 And $K(G)\# K[G]$ is the braided product $\mathscr{H} = K(G) \propto K[G]$,
 where $K(G)$ and $K[G]$ are braided commutative Yetter-Drinfel'd $K(G)\bowtie K[G]$-module algebras by restriction, i.e.,
 the action is given by
 \begin{eqnarray*}
 && (\delta_{p} \bowtie q)\cdot \delta_{p'} = \delta_{qp'^{-1}p} \delta_{p'q^{-1}}, \\
 && (\delta_{p} \bowtie q)\cdot q' = \delta_{p^{-1}, qq'q^{-1}} qq'q^{-1},
 \end{eqnarray*}
 and coaction $\rho$: $A \rightarrow M_{0}(\mathscr{D} \otimes A)$ and $B \rightarrow M_{0}(\mathscr{D} \otimes B)$
 is given by $\rho(\delta_{p})=\sum_{t\in G} \delta_{t^{-1}p} \bowtie e \otimes \delta_{t}$,
 $\rho(q)=\sum_{t\in G} \delta_{t} \bowtie q \otimes q$,
 for$p, p', q, q' \in G$ and $\delta_{p}, \delta_{p'} \in K(G)$.

\section{Some special cases}
\def\theequation{\thesection.\arabic{equation}}
\setcounter{equation}{0}

 In this section, we apply our results as above to the usual Hopf algebras (i.e., multiplier Hopf algebra has an identity), and derive some interesting results.
 \\

 Let $B$ be a (infinite dimensional) co-Frobenius Hopf algebra with a left integral $\varphi$,
 and $A$ be the dual multiplier Hopf algebra shown in \cite{Zh99}.
 Then by \cite{DV04} or Corollary 3.6 in \cite{YW11a}, let $\Phi_{\alpha}=\iota=\Phi_{\beta}$,
 we can get the Drinfel'd double $\mathscr{D}=A\bowtie B$ with structures
 given by the following formulas:
 \begin{eqnarray*}
 &&(a\bowtie b)(a'\bowtie b') = a(b_{(1)} \blacktriangleright a' \blacktriangleleft S^{-1}(b_{(3)})) \bowtie b_{(2)}b',\\
 && \Delta(a\bowtie b) = \Delta^{cop}(a)(b_{(1)} \otimes b_{(2)}), \\
 && \varepsilon(a\bowtie b) = \varepsilon_{A}(a)\varepsilon_{B}(b),\\
 && S(a\bowtie b) = S(b_{(3)}) \blacktriangleright S^{-1}(a) \blacktriangleleft b_{(1)} \otimes S(b_{(2)})
 \end{eqnarray*}
 for any $a\in A$, $b\in B$.
 \\

 Let the Heisenberg double $\mathscr{H}=A\# B$ with the multiplication as (\ref{3.1}) and endow $\mathscr{H}$ with
 the $\mathscr{D}$-module and comodule structures as (\ref{3.2}) and (\ref{3.3}). Then by Theorem 3.4 and Corollary 3.8,
 we get the following result in the form of a theorem, which gives an answer to the question introduced in the introduction.
 \\

 \textbf{Theorem \thesection.1}
 Let $B$ be a co-Frobenius Hopf algebra with a left integral $\varphi$, and $A$ be the dual regular multiplier Hopf algebra.
 Then $\mathscr{H}$ is a (left-left) $\mathscr{D}$-Yetter-Drinfel'd module algebra,
 and moreover $\mathscr{H}$ is a braided $\mathscr{D}$-commutative algebra.
 \\

 \textbf{Example \thesection.2}
 Let $C$ be an infinite cyclic group with generator $c$ and let $m$ be a positive integer.
 Let $i\in \mathbb{N}$, the set of natural integers  and $\lambda \in \mathbb{C}$ such that $\lambda^{i}$
 is a primitive $m$th root of $1$.
 Then we recall from \cite{YW11a} that the Hopf algebra $B$ is the algebra with generators $c$ and
 $X$ satisfying relations: $cX=\lambda Xc$ and $X^m=0$.
 The Hopf algebra structure on $B$ is given by
 \begin{eqnarray*}
 && \Delta(c)=c \otimes c, \qquad \Delta(X)=c^i \otimes X + X \otimes 1, \\
 && \varepsilon(c)=1, \qquad  \qquad \varepsilon(X)=0, \\
 && S(c)=c^{-1}, \qquad \quad S(X)=-c^{-i}X.
 \end{eqnarray*}

 In [\cite{DV04}, 2.2.1], the authors construct the multiplier Hopf algebra $A=\widehat{B}$ with
 the linear basis $\{\omega_{p,0}Y^{l}\mid  p\in \mathbb{Z}, l\in \mathbb{N}, l<m \}$.
 The multiplication  and the comultiplication are defined so that $\< A, B\> $ is a multiplier Hopf algebra pairing.
 For the details, the product in $A$ is given by the formula
 $\omega_{p, q}\omega_{k, l}=\delta_{p-k, il}({}^{l+q}_{q})_{\lambda^{-i}}\omega_{k, l+q}$
 and the multiplier Hopf structure of $B$ is given by
 \begin{eqnarray*}
 && \Delta(\omega_{p, 0})=\sum_{k\in \mathbb{Z}} \omega_{k, 0} \otimes \omega_{p-k, 0},
  \quad \Delta(Y)=D\otimes Y+Y \otimes 1, \\
 && \varepsilon(\omega_{p, 0})=\delta_{p, 0}, \qquad \qquad \quad \varepsilon(Y)=0, \\
 && S(\omega_{p,0})=\omega_{-p, 0},  \qquad  \qquad  S(Y)=-D^{-1}Y,
 \end{eqnarray*}
 where $D=\sum_{j\in \mathbb{Z}}\lambda^{j}\omega_{j, 0}$ and
 $Y=\sum_{s\in \mathbb{Z}}\lambda^{s}\omega_{s,1}$. Notice that
 $DY=\lambda YD$, $Y^{m}=0$ and $D\omega_{k, 0}=\lambda^{k}\omega_{k, 0}=\omega_{k, 0}D$.

 Define Heisenberg double $\mathscr{H}=A\# B$ as follows,
 \begin{eqnarray*}
 && (\varepsilon \# c^{i})(\omega_{p,0} \# 1)=\omega_{p-i,0} \# c^{i}, \quad (\varepsilon \# c^{i})(Y \# 1)=Y\# c^{i}, \\
 && (\varepsilon \# X)(\omega_{p,0} \# 1)=\omega_{p-i,0} \# X, \quad  (\varepsilon \# X)(Y \# 1)=Y\# X+ D\# 1,
 \end{eqnarray*}
 and Drinfel'd double $\mathscr{D}=A\bowtie B$ as Example 2.8 in \cite{YW11a}(only need to let $\alpha=\beta=\gamma=\delta=\iota$),
 then we can get $\mathscr{H}$ is a $\mathscr{D}$-Yetter-Drinfel'd module algebra,
 and it is braided commutative.
 \\

 In the follwing, we consider the case that all the two regular multiplier Hopf algebras $A$ and $B$ have identities $1_{A}$ and $1_{B}$ respectively, i.e. $A$ and $B$ are Hopf algebras by Theorem 2.4 in \cite{V98}.

 Let $S_{A}$ and $S_{B}$ be the be bijective antipodes of Hopf algebras $A$ and $B$ respectively,
 and $(A, B, \<\cdot, \cdot\>)$ be a Hopf dual pairing, then we can get Drinfel'd double $\mathscr{D}=A\bowtie B$
 with structures
 \begin{eqnarray*}
 && (a\bowtie b)(a'\bowtie b') = \<a'_{(1)}, S^{-1}_{B}(b_{(3)})\> (aa'_{(2)}\otimes b_{(2)}b') \<a'_{(3)}, b_{(1)}\>, \\
 && \Delta(a\bowtie b) = (a_{(2)}\bowtie b_{(1)}) \otimes  (a_{(1)}\bowtie b_{(2)}), \\
 && \varepsilon(a\otimes b) = \varepsilon_{A}(a)\varepsilon_{B}(b), \quad \mbox{and}\\
 && S(a\otimes b)= \<a_{(1)},  b_{(3)}\> (S^{-1}_{A}(a_{(2)})\otimes S_{B} (b_{(2)})) \< S^{-1}_{A} (a_{(3)}), b_{(1)} \>
 \end{eqnarray*}
 and Heisenberg double $\mathscr{H}=A\# B$ with multiplication
 \begin{eqnarray*}
 && (a \# b)(a' \# b') = \< a'_{(2)}, b_{(1)}\> aa'_{(1)} \# b_{(2)}b'.
 \end{eqnarray*}
 for any $a, a' \in A$ and $b, b' \in B$.

 Define $\mathscr{D}$-module action on $\mathscr{H}$
  \begin{eqnarray*}
 && (a \bowtie b)\cdot (a' \# b') = a_{(3)} (b_{(1)} \blacktriangleright a') S^{-1}(a_{(2)}) \#
 (b_{(2)} b' S(b_{(3)})) \blacktriangleleft S^{-1}(a_{(1)}).
 \end{eqnarray*}
 and comodule action $\rho: \mathscr{H} \longrightarrow \mathscr{D} \otimes \mathscr{H}$,
 $\rho(a\# b)=a_{(2)}\bowtie b_{(1)} \otimes a_{(1)}\# b_{(2)}$.
 Then we can get
 \\

 \textbf{Corollary \thesection.3}
 $\mathscr{H}=A\# B$ is a (left-left) $\mathscr{D}=A\bowtie B$-Yetter-Drinfel'd module algebra,
 and $\mathscr{H}$ is a braided $\mathscr{D}$-commutative algebra.
 \\

 Furthermore, if $B$ is a finite dimensional Hopf algebra, then the antipode is bijective, and we can construct its duality $B^{*}$,
 which is also a Hopf algebra satisfying the condition of a pairing. So we can get a corollary, which is the main results in \cite{S11}.
 \\

 \textbf{Corollary \thesection.4}
 $\mathscr{H}(B^{*})=B^{*}\# B$ is a (left-left) $\mathscr{D}(B)=B^{*}\bowtie B$-Yetter-Drinfel'd module algebra,
 and $\mathscr{H}$ is a braided $\mathscr{D}$-commutative algebra in terms of the braiding of Yetter-Drinfel'd module.

\section*{Acknowledgements}

 The authors would like to thank the referee for his/her valuable comments.
 The work was partially supported by the NNSF of China (No. 11226070, 11571173), the NJAUF (No. LXY201201019, LXYQ201201103)
 and NSF for Colleges and Universities in Jiangsu Province (No. 11KJB110004).

\end {document}